# A “Travelling Salesman’s Paths” within $n \times n$ ($n$ = 3, 4, 5) Magic Squares

Peyman Fahimi,[1*] Walter Trump,[2] Chérif F. Matta,[3] Alireza Ahmadi Baneh[4]

[1] *Department of Mathematics and Statistics, Dalhousie University, Halifax, NS B3H4R2, Canada.*
** Tel.: +1-(438)-229-2038, E-mail: fahimi@dal.ca.*
[2] *Department of Physics, Gymnasium Stein, 90547 Stein, Germany.*
[3] *Department of Chemistry and Physics, Mount Saint Vincent University, Halifax, NS B3M2J6, Canada.*
[4] *Department of Applied Mathematics, Tarbiat Modares University, Tehran, Tehran, Iran.*

**Abstract**
Intriguing symmetries are uncovered regarding all magic squares of orders 3, 4, and 5, with 1, 880, and 275,305,224 distinct configurations, respectively. In analogy with the travelling salesman problem, the distributions of the total topological distances of the paths travelled by passing through all the vertices (matrix elements) only once and spanning all elements of the matrix are analyzed. Symmetries are found to characterise the distributions of the total topological distances in these instances. These results raise open questions about the symmetries found in higher-order magic squares and the formulation of their minimum and maximum total path lengths.

## 1. Introduction

Magic squares have captivated mathematicians, historians, and enthusiasts for centuries with their intricate patterns and inherent mathematical properties [1–7]. The concept of magic squares can be traced back thousands of years to ancient civilizations such as China [8], India [9], and the Middle East [9]. Over the centuries, magic squares have captured the imagination of mathematicians including, for example, the Persian mathematician Abu al-Wafa'al-Buzjani [9] up to Leonard Euler [10]. These intriguing structures, are square matrices filled with distinct integers such that they possess a remarkable property, and that is that the sum of any row or column as well as its trace are all equal.

A classic or natural magic square of order $n$ is an array of distinct positive integers ranging from 1 to $n^2$. This arrangement ensures that the sum of any $n$ numbers in a horizontal, vertical, or main diagonal line consistently equals a specific value known as the magic constant. The mathematical expression for the magic constant is $M_n = \frac{n}{2}(n^2 + 1)$. We define two subgroups of magic squares that will be explored later in this paper. An $n \times n$ classic magic square earns the label "associative" when each pair of numbers symmetrically opposite the center sums to $n^2 + 1$ [4]. Moreover, a classic magic square earns the designation of a pandiagonal magic square if the sum of all its diagonals, including those derived by "wrapping around" the edges, equals the same magic constant [4].

The number of possible (classic) magic squares of a given dimension explodes very quickly with matrix size [11–14]. Let’s define an “irreducible” square as one with the smallest denominator of the infinite number of corresponding squares obtained by multiplication with a constant. A 3×3 magic square has only one possible irreducible arrangement, there are 880 irreducible 4×4 magic squares, there are 275,305,224 irreducible magic square of size 5×5, and an estimated *ca.* $2.4 \times 10^{110}$

of size 10×10 [11]. The law determining the number of magic squares of a given type given the matrix size remains unknown and an active area of research [11].

Magic squares, besides their mathematical elegance, aesthetic appeal, and remarkable properties, may have some applications as well. Magic squares have been employed in classical mechanics [15,16], electrostatics [17–20], and even quantum mechanics [21,22].

Fahimi [17] recently proposed a method to make the 4×4 magic squares, "float" in the air like magic. The trick is to think of the numbers in the arrangement as if they were electric charges on a 2D lattice. They use a computer simulation with a smart controller to show how these imaginary charges behave in a symmetrical way when placed in an electric field.

The representation of magic squares in binary format reveals underlying patterns [23], which may hold potential implications in binary systems within the realm of physics, such as the 2D Ising model [24]. In a recent study by Fahimi [23], the author introduces a method for exploring similarities within magic squares by assigning two distinct colors to the even and odd numbers. This technique encompasses various analyses such as rotation/reflection, Principal Component Analysis (PCA), and Linear Discriminant Analysis (LDA). Initially, the magic squares undergo conversion into their binary format. The subsequent step involves flattening each binary matrix into an $n^2$-dimensional vector, where '$n$' denotes the order of the square. Following this transformation, a dimensionality reduction technique, such as PCA, is applied to condense the information into two components, offering a more accessible representation. However, a rigorous mathematical and analytical exploration is necessary as a follow-up to this numerical investigation. This entails tasks like computing covariance matrices and performing eigenvalue decomposition for each binary magic square. Such a comprehensive approach unravels valuable insights into the underlying principles, enriching our comprehension of the intricate dynamics at play in the realm of magic squares.

In the present study, matrix elements within a magic square are taken as a numerical label of various cities. The objective is to investigate the symmetries in the trajectories of a traveler moving from city 1 to city $n^2$, where $n$ represents the order of the magic square (with $n$ = 3, 4, and 5 being examined in this paper). By examining the "traveling salesman" trajectories from a city to the next, the aim is to uncover and analyze the inherent symmetrical patterns present in these magical arrangements.

## 2. Methodology

The methodology, implemented in Python, utilizes data handling libraries such as Pandas and NumPy. A portion of the code related to 4×4 magic squares is included in the Appendix, while other codes, including the code for 5×5 magic squares, as well as the results and plots, are included in the supplementary material.

The code begins by loading magic squares from an Excel file and defining coordinates for cities within the square grid. The loop iterates through each row of the DataFrame containing the magic squares. For each iteration, the code extracts a magic square from the DataFrame and converts it into a list. It then determines the maximum number present in that particular magic square. The next step involves creating a dictionary ("city_mapping") that associates each number in the magic square with its corresponding city coordinates. The cities dictionary holds the predefined coordinates for each city. Using the created "city_mapping" dictionary, the code generates a path for the traveler by mapping the numbers in the magic square to their respective city coordinates. The resulting "*x*" and "*y*" lists contain the coordinates of the traveler's path.

For each magic square, the trajectory of the traveler is plotted, showcasing the movement between cities with arrows. The distances between consecutive cities are calculated, utilizing the Euclidean distance formula. This distance calculation serves as a fundamental metric for analyzing spatial relationships within the magic square arrangement. The resulting plots visually represent both the trajectory and the calculated distances, offering valuable insights into the symmetrical dynamics embedded in the mathematical constructs of magic squares.

## 3. Distance Measurement and Symmetry Analysis of City Trajectories

Let us consider a 3×3 magic square. As depicted in **Fig. 1 (left)**, we envision a distribution of 9 cities on a 3×3 square lattice whereby each neighboring pair of cities is separated by a unit horizontal or vertical distance. The arrangement of the cities precisely mirrors the arrangement of the numbers within a 3$^{rd}$ order magic square, resulting in a magical path for the "travelling salesman".

All possible inter-pair distances are analyzed for their underlying pattern. For instance, the distance between city 1 and city 2 is $\sqrt{(1^2+2^2)}=\sqrt{5}$, while the distance between city 4 and city 5 is $\sqrt{2}$, and so on. In the particular path depicted in **Fig. 1 (right)**, the distances pattern exhibits *reflective (σ) symmetry* about a mirror plane intersection the plane of the figure vertically between abscissa steps 3 and 4. The total distance of the sum of all paths is approximately 13.77. It is worth clarifying at this point that the paths considered in this study must align with the *order of occurrence of consecutive magic numbers* and not any arbitrary path. In other words, we adhere to the constraint that moving from, for instance, city 1 to city 5 is not permissible. One might, of course, question the possibility of exploring alternative routes, such as moving from 1 to 5 and then to 9, but such deviations are outside the scope of this initial study and may be the subject of a future investigation.

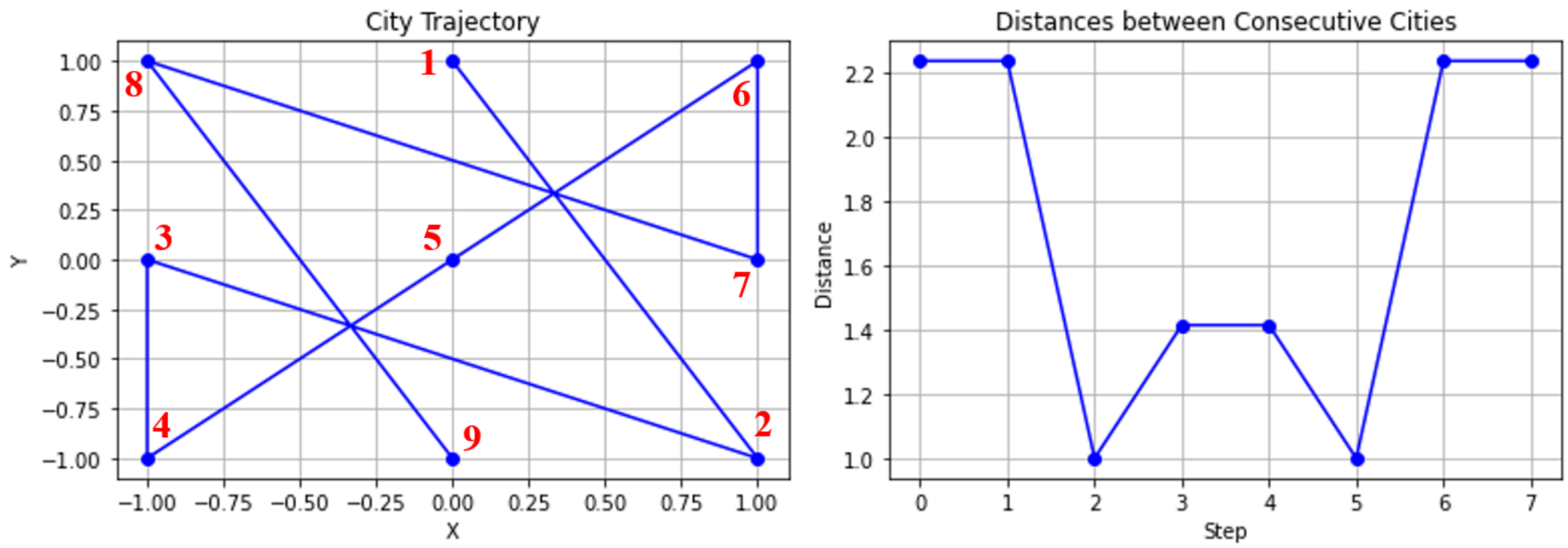

**Figure 1:** 3×3 magic travelling salesman trajectory (the plot on the right shows the σ-(reflective) symmetry on an imaginary mirror plane perpendicular to the vertical line midway between step 3 and 4. The reader should note that the initial step is labeled as step zero denoting the distance from city 1 to city 2, which is equivalent to $5^{½}$ or ~ 2.2.

The complexity of paths in 4×4 magic squares represent a greater challenge in discerning their underlying symmetries in terms of topological distance patterns. Among the 880 order 4 magic squares, it is observed that the maximum and minimum values of the total path distances are approximately 42.76 and 20.31, respectively, while the average total path distance over all 880 magic squares is *ca*. 33.94.

**Fig. 2a** and **2b** illustrate the magic squares with the minimum and maximum total distances, respectively, while **Fig. 2c** displays the distribution of total distances across all 880 magic squares of order 4. **Fig. 2a** represents the shortest pathway for the traveler, while **Fig. 2b** depicts the longest pathway among all 880 routes. Furthermore, the analysis shows that the shortest total trajectory is 20.3, with an average distance of 1.35 between each pair of consecutive cities (as illustrated in Fig. 2a and the supplementary material). In contrast, the longest trajectory extends to 42.7, with an average distance of 2.85 between each pair of consecutive cities (as depicted in Fig. 2b and the supplementary material). The distribution of total distances across all 880 magic squares of order 4 is presented in Fig. 2c.

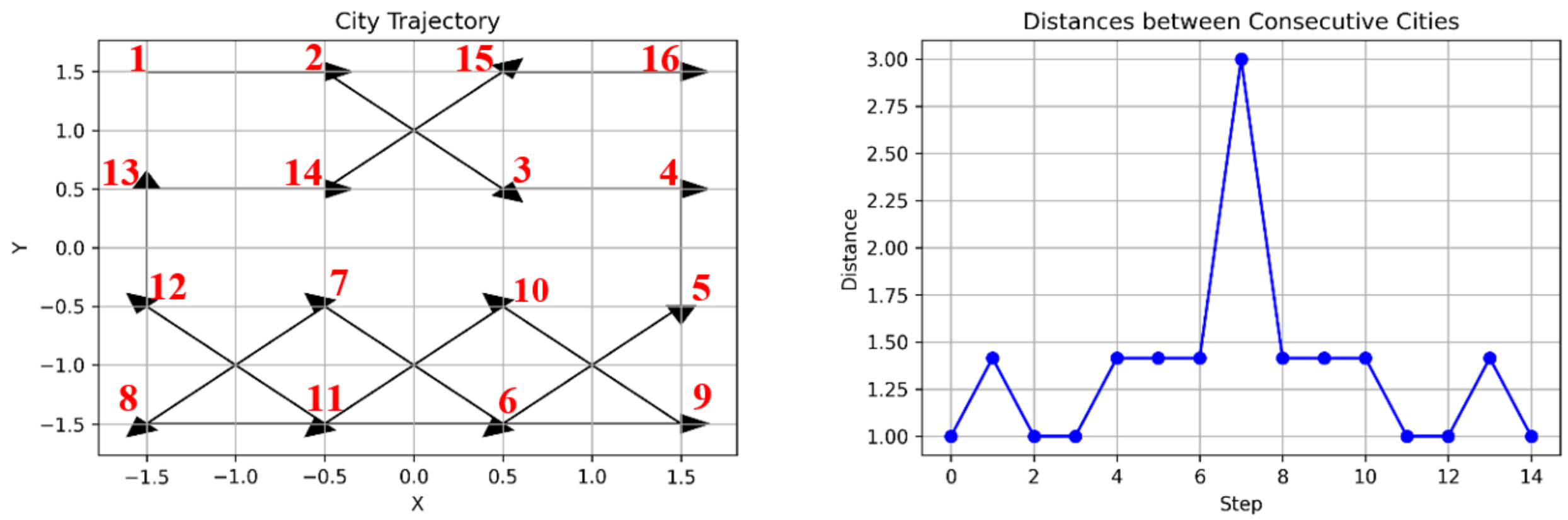

(a)

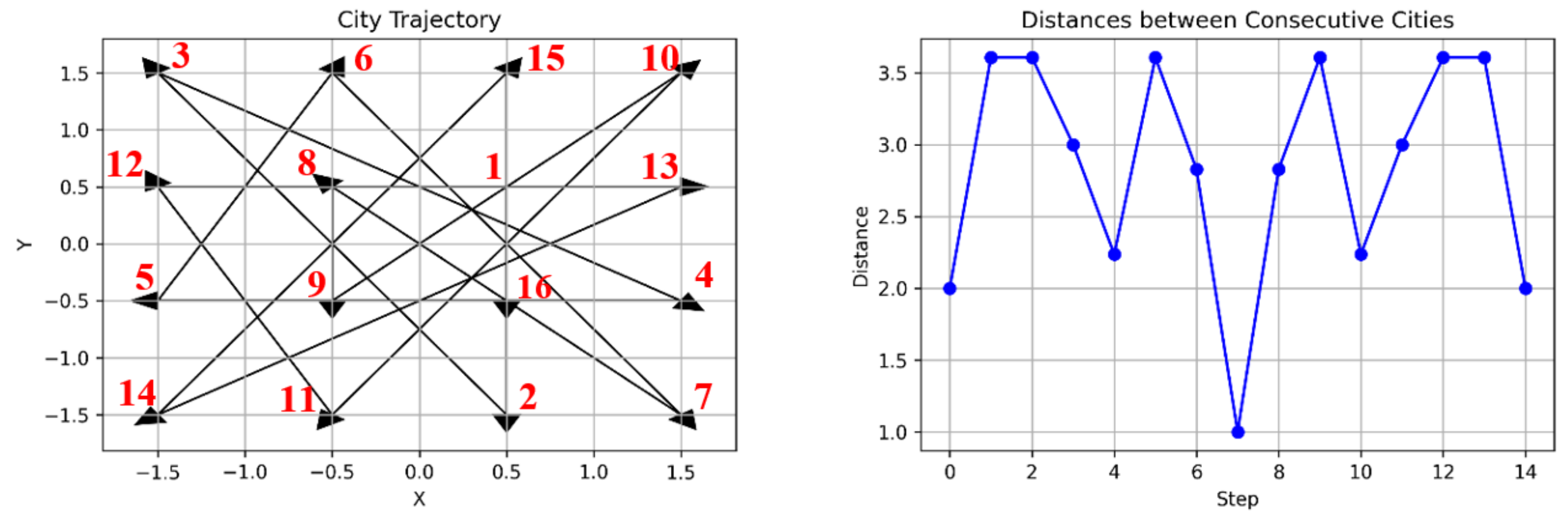

(b)

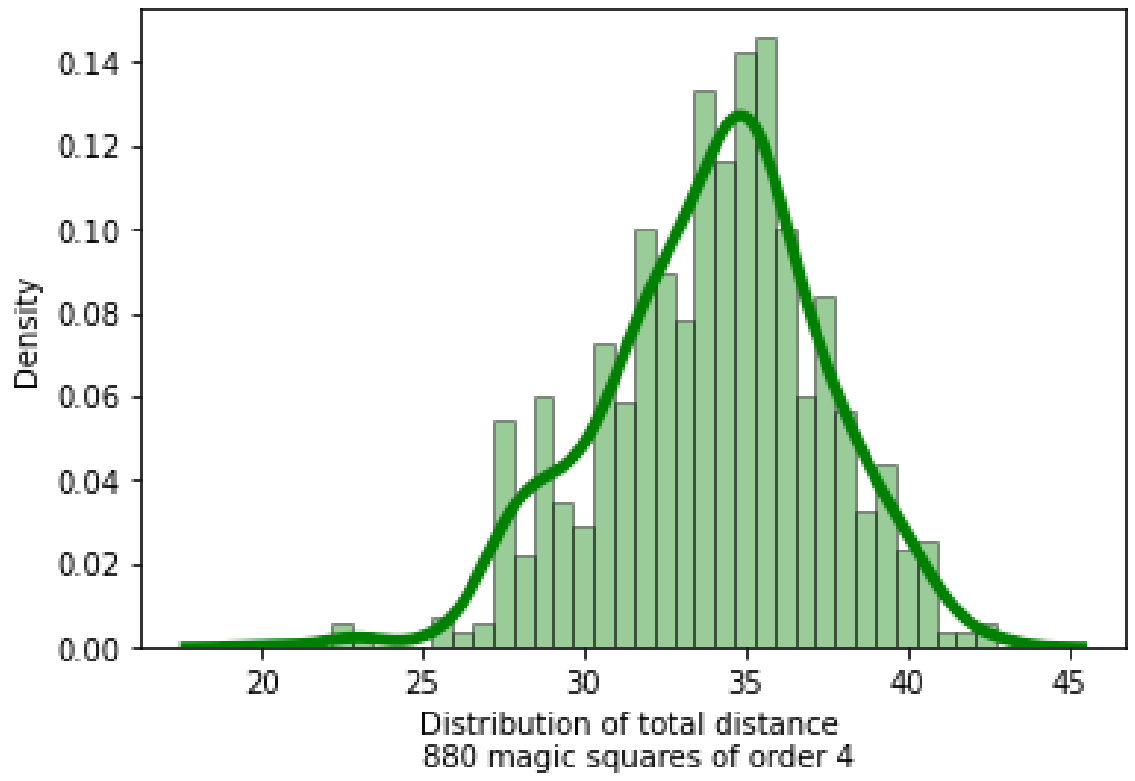

(c)

**Figure 2:** Panel (a) illustrates the magic square with the shortest total trajectory, panel (b) showcases the magic square with the longest total trajectory, and panel (c) displays the histogram depicting the distribution of total distances across all 880 magic squares of order 4.

Within the set of 880 magic squares of order 4, 414 display reflective geometric symmetry in their distance patterns. This symmetry is observed around a mirror plane intersecting the figure's plane vertically at step 7. In the pursuit of understanding the origins of such symmetries, an in-depth exploration of the 880 magic squares with dimensions 4×4 was conducted, leveraging a classification system devised by H. E. Dudeney [3,25]. These squares were systematically organized into 12 distinct groups based on the patterns formed by the 8 complement pairs, as illustrated in the left panel of Figure 3 (refer to **Fig. 3a** and **b**). Here, a complement pair denotes two numbers whose sum equals $n^2 + 1$, where $n$ represents the order of the magic square; in the case of order 4, this sum is 17.

The classification according to Dudeney's scheme provides a structured framework for deciphering the underlying patterns that contribute to the emergence of these symmetrical arrangements [3,25]. It is found that all magic squares belonging to group 3, which consists of 48 associative magic squares with symmetrical number placement around the center point, as well as group 6, comprising 304 semi-pandiagonal and simple magic squares with symmetrical number placement across the center line, demonstrate reflexive symmetry in their distance patterns. **Fig. 3** showcases a selection of examples from these groups. A comprehensive collection of all examples is available as Supplementary Information. Out of the 880 distance patterns, it is important to note that several are redundant. Instead, there are 112 emerging patterns that exhibit repetition, resulting in a total of 768 unique distance patterns. The reader is reminded that each magic square corresponds to a unique trajectory and total distance pattern, given the constraint of having to follow the sequential numbering of the cities. Thus, we have, for example, a total of 880 trajectories for the 4×4 magic squares within which 112 are identical despite being associated distinct magic squares.

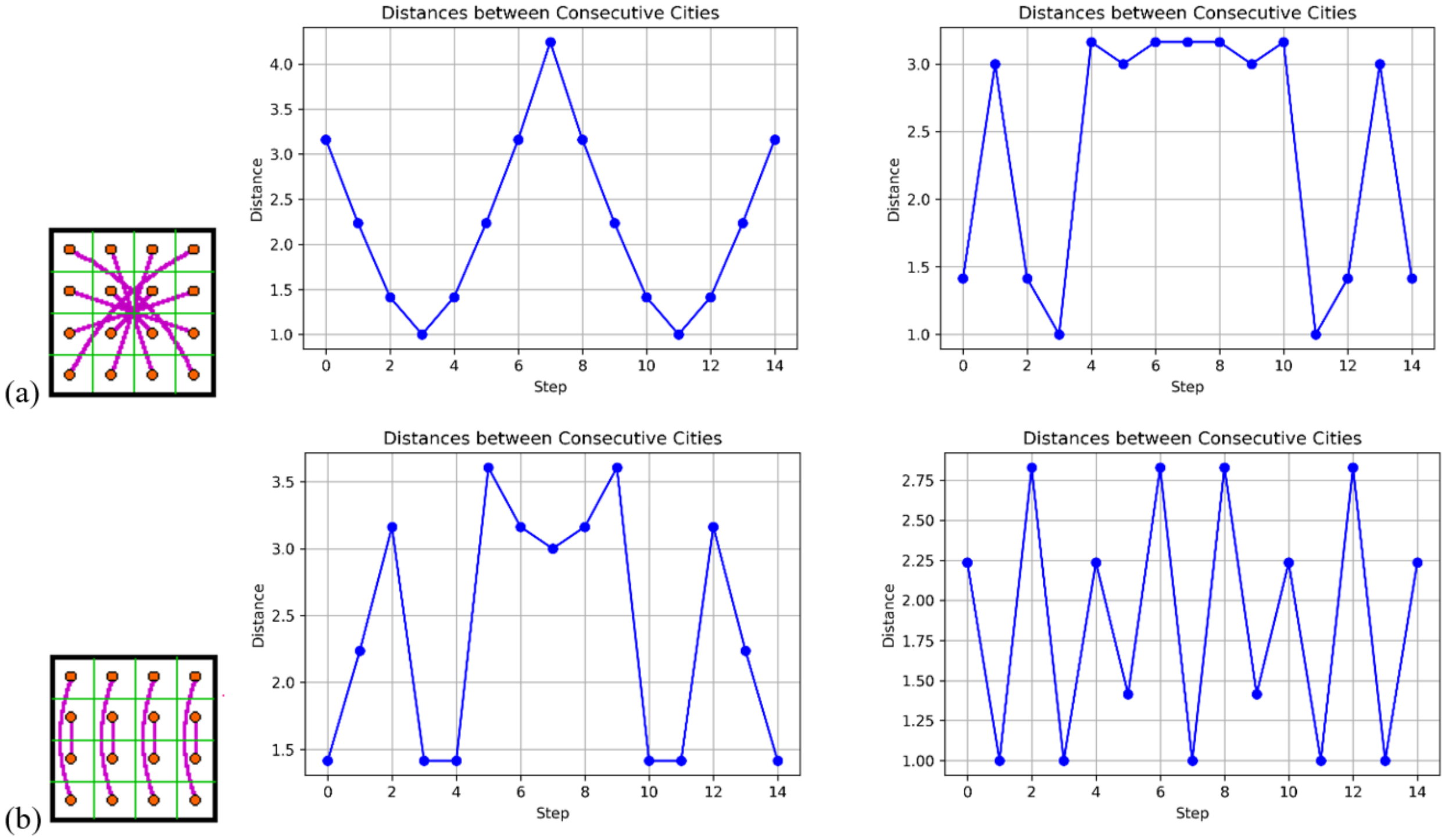

**Figure 3:** Illustration of symmetrical distance patterns in the traveler's path, showcasing examples from Dudeney group 3 (panel a) and group 6 (panel b). The purple plots (left plot in each panel) were obtained from the website of Harvey Heinz *[25]*. The copyright permission has been obtained from the original source, the book "Amusements in Mathematics" by Henry Ernest Dudeney [3], which is published by Dover Publications.

The remaining distance patterns, comprising 466 cases (= 880 – 414), exhibit various types of symmetries, including local symmetry, periodicity (translational symmetry), and partial symmetry. Local symmetry (which 252 patterns exhibit this type of symmetry) refers to the presence of symmetric patterns or characteristics within a specific region or subset of the distance pattern. It implies that certain portions of the pattern exhibit mirror-like or rotational symmetry within themselves. Periodicity or translational symmetry refers to the recurrence or repetition of certain patterns or motifs at regular intervals within the distance pattern. It indicates the existence of a periodic structure or arrangement in the pattern. Partial symmetry suggests that the distance pattern possesses symmetrical elements or features, but not all aspects of the pattern exhibit symmetry. It indicates the presence of symmetry in some parts or components of the pattern while other parts may lack symmetry. In the majority of cases, achieving reflexive symmetry, local symmetry, or periodicity can be accomplished by repositioning a few points (1 to 3 points) in the distance pattern of partial symmetry. These classifications are primarily based on qualitative assessment and may exhibit overlap in certain cases. **Fig. 4** depicts a few examples of different types which are listed in their entirety in the Supplementary Information.

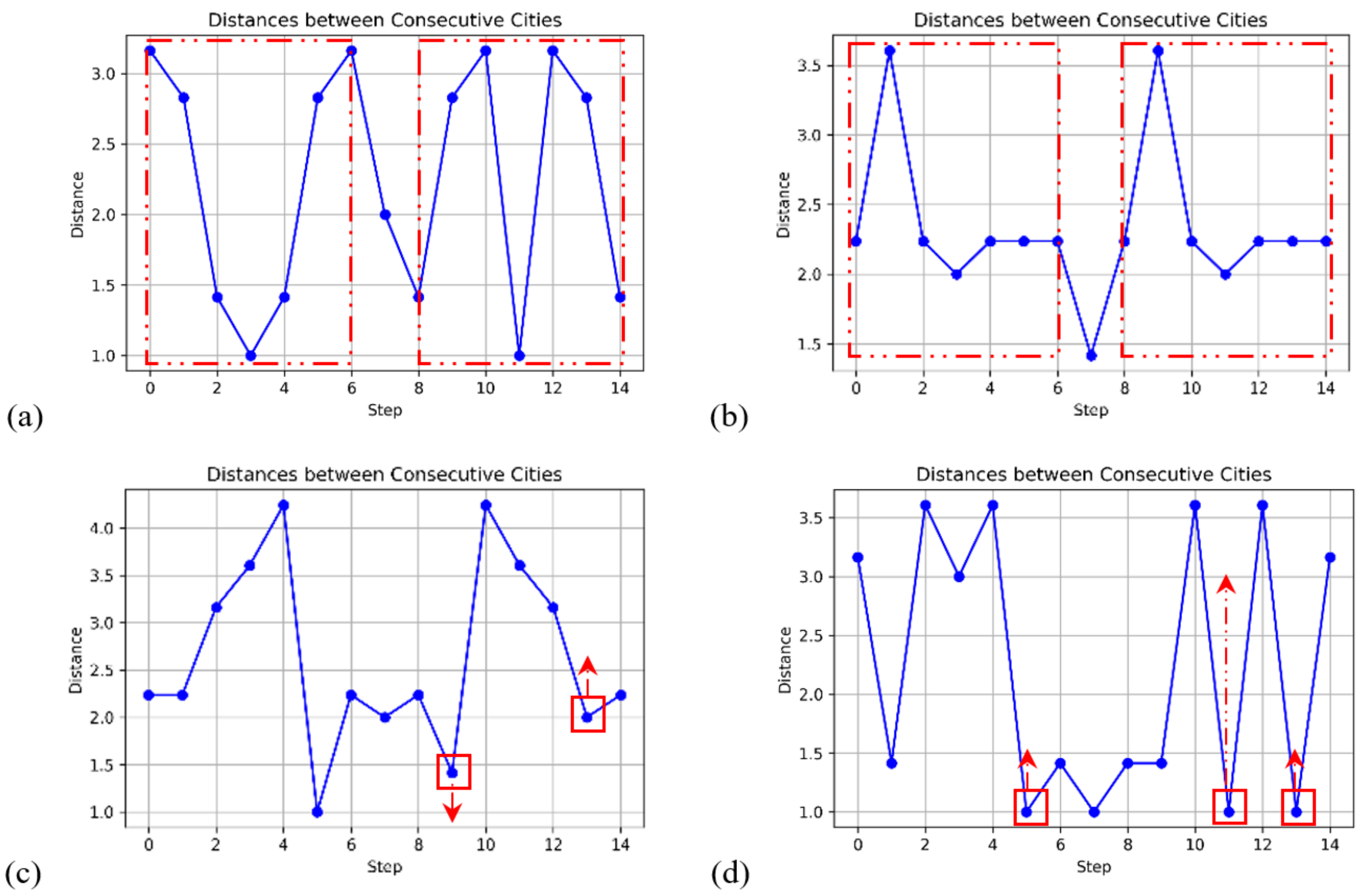

**Figure 4:** Examples of distance patterns: local symmetry (a), periodicity (b), partial symmetry (c, d).

The presence of partial symmetry in distance patterns can be attributed to the specific arrangement of the cities within the square. The topological distances of different pairs of cities introduce variations that disrupt the symmetry embodied in the magic square by construction. The arrangement of cities within the magic square creates complex interdependencies, leading to partial symmetric patterns in the distances between certain pairs of cities. These patterns may arise due to the specific positioning of cities along diagonal lines, within certain quadrants of the square, or other geometric relationships. It is crucial to emphasize that while the city numbers demonstrate symmetry by adhering to the arrangement of numbers within a magic square, the actual distances between cities are governed by geometric considerations and may not necessarily align with the same level of symmetry. This creates an interesting interplay between the overall symmetric structure of the magic square and the partial symmetries observed in the distance patterns. By analyzing and understanding these partial symmetries, further insights into the underlying principles and properties of magic squares and their associated distance patterns can be gained.

Investigating the trajectory patterns and distance symmetries within 5×5 magic squares poses significant computational challenges due to the vast number of distinct magic squares of order 5, which amounts to 275,305,224 [11]. We calculated the total path length and their respective frequencies, of all these squares, as illustrated in **Fig. 5**. The resulting pattern resembles a Gaussian-type distribution, with the maximum frequency of 32,824 occurring at a total path length of approximately 65.89. The minimum, average, median, maximum, and standard deviation of the total path length are 37.47, 65.49, 65.58, 88.04, and 4.64, respectively. Note that, the diagrams in Fig. 5 and Fig. 7 do not display the exact total path length *S* and their frequencies. Rather, they show the number of squares with total path lengths that fall within a small interval around *S*, specifically [*S* - 0.005, *S* + 0.005].

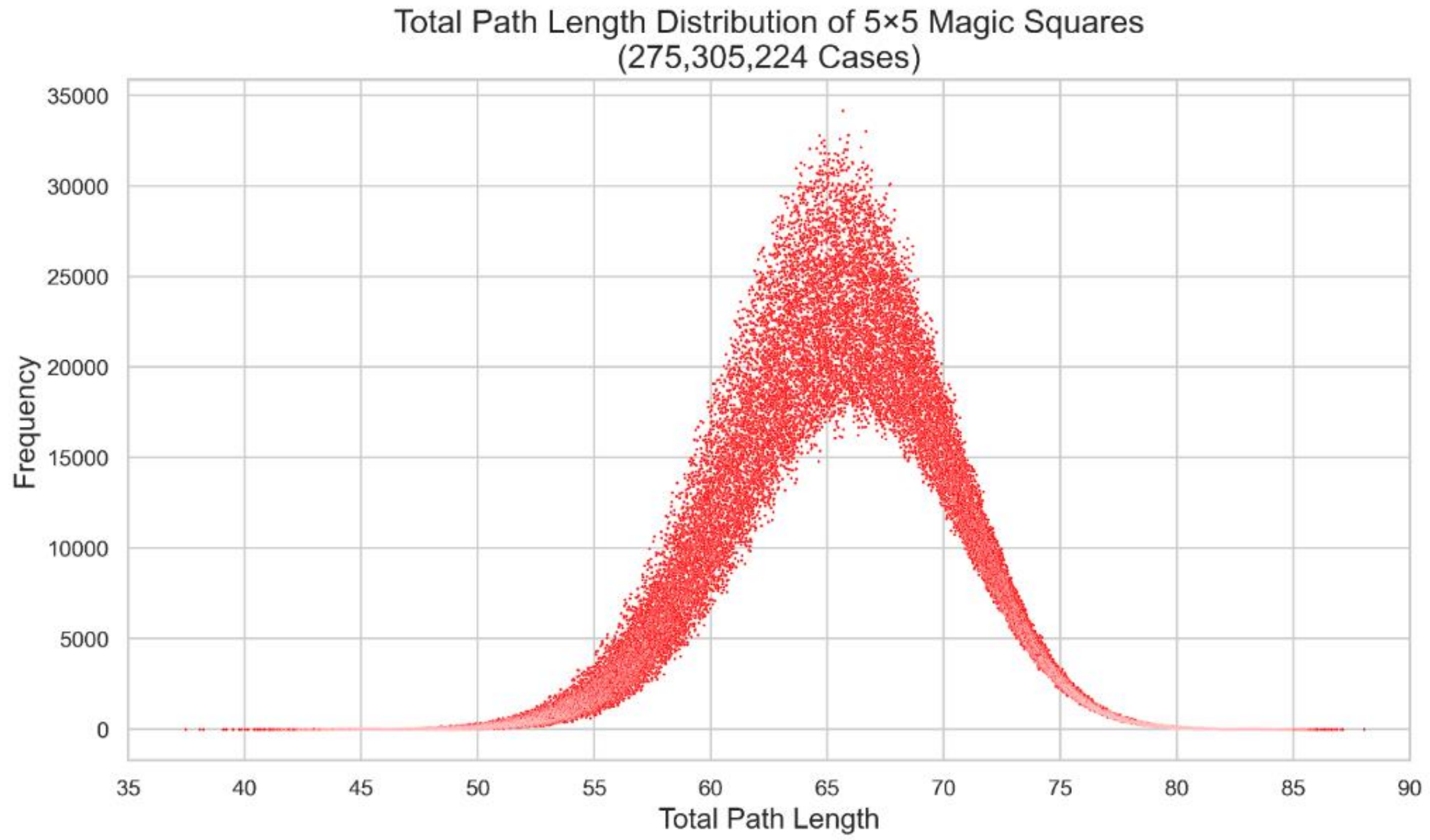

**Figure 5:** Gaussian-type distribution of total path lengths and their corresponding frequencies for the traveling salesman in 5x5 magic squares.

Our analysis demonstrates that the distance pattern of the traveler's path in all 48,544 instances of associative magic squares of order 5 exhibits a remarkable reflexive geometrical symmetry. Additionally, there are 5,848 non-associative 5×5 magic squares exhibiting reflexive path symmetry. This group of squares also features the number 13 at the center, similar to the associative magic squares. **Fig. 6** illustrates some notable examples of these patterns.

In a 5×5 magic square, any number from 1 to 25 can occupy the center. Therefore, we categorized all magic squares of order 5 into 25 groups, each defined by a specific number in the center. The results for having *x* in the center are identical to those with (26 - *x*) in the center. Consequently, we plotted the total path length distribution for the first 13 groups, with center numbers ranging from 1 to 13. **Fig. 7** illustrates these distributions, with the maximum frequency and average total path length provided in the legend of each plot. As shown, the average total path length across all groups hovers around 64. Notably, in most cases, the smaller the center number, the lower the height of the Gaussian distribution (i.e., the maximum frequency), though there are exceptions.

Assuming the traveling salesman seeks the minimum total path length (MTPL) for each order of magic squares, the following outlines the steps of the shortest trajectory obtained from numerical analysis of all magic squares of orders 3, 4, and 5:

for *n* = 3 (only possible path):
TPL = $( \mathbf{2} \times \sqrt{(1^2)} ) + ( 2 \times \sqrt{(1^2+1^2)} ) + ( 4 \times \sqrt{(1^2+2^2)} ).$

for *n* = 4 (minimum total path length):
MTPL = $( \mathbf{6} \times \sqrt{(1^2)} ) + ( 8 \times \sqrt{(1^2+1^2)} ) + ( 1 \times \sqrt{(3^2)} ),$

for $n = 5$ (minimum total path length):

$$MTPL = ( \mathbf{10} \times \sqrt{(1^2)} ) + ( 8 \times \sqrt{(1^2+1^2)} ) + ( 4 \times \sqrt{(1^2+2^2)} ) + ( 2 \times \sqrt{(2^2+3^2)} ).$$

Even with the results of the shortest trajectories for magic squares of orders 3, 4, and 5, it remains challenging to propose a general formula for predicting the minimum total path length for higher-order magic squares. Additionally, calculating the total path length (TPL) for higher-order squares is computationally intensive; for instance, there are 17,753,889,197,660,635,632 possible magic squares of order 6 [11]. Therefore, developing a mathematical formula to calculate the minimum total path length (MTPL) for any order remains an open problem for future research.

| | | | | |
|---|---|---|---|---|
| 5 | 22 | 14 | 18 | 6 |
| 17 | 11 | 2 | 25 | 10 |
| 7 | 23 | 13 | 3 | 19 |
| 16 | 1 | 24 | 15 | 9 |
| 20 | 8 | 12 | 4 | 21 |

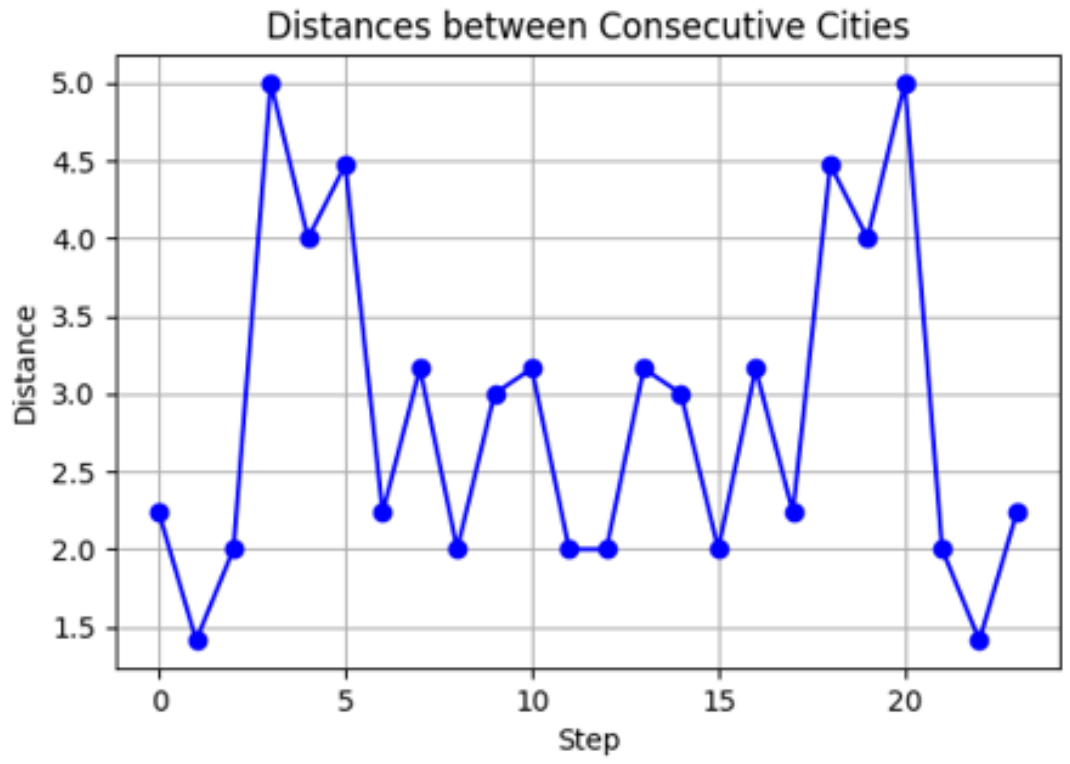

(a)

| | | | | |
|---|---|---|---|---|
| 8 | 23 | 15 | 7 | 12 |
| 16 | 4 | 1 | 24 | 20 |
| 21 | 17 | 13 | 9 | 5 |
| 6 | 2 | 25 | 22 | 10 |
| 14 | 19 | 11 | 3 | 18 |

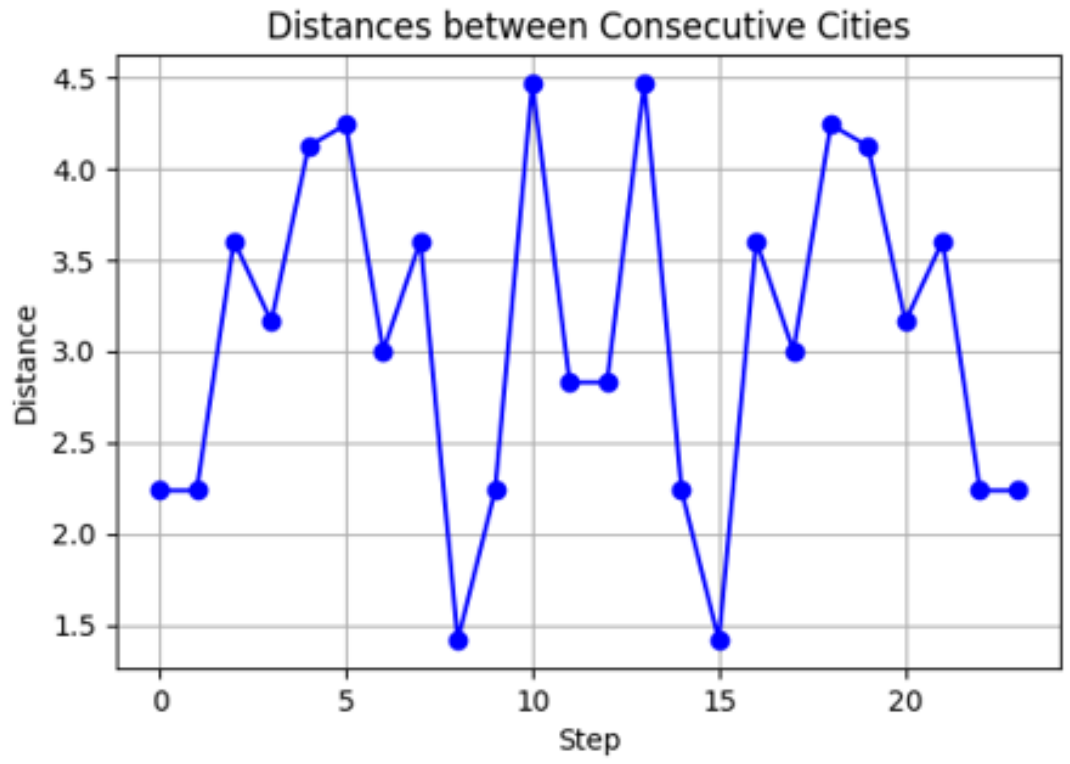

(b)

| | | | | |
|---|---|---|---|---|
| 7 | 10 | 20 | 17 | 11 |
| 3 | 21 | 12 | 25 | 4 |
| 18 | 24 | 13 | 2 | 8 |
| 22 | 1 | 14 | 5 | 23 |
| 15 | 9 | 6 | 16 | 19 |

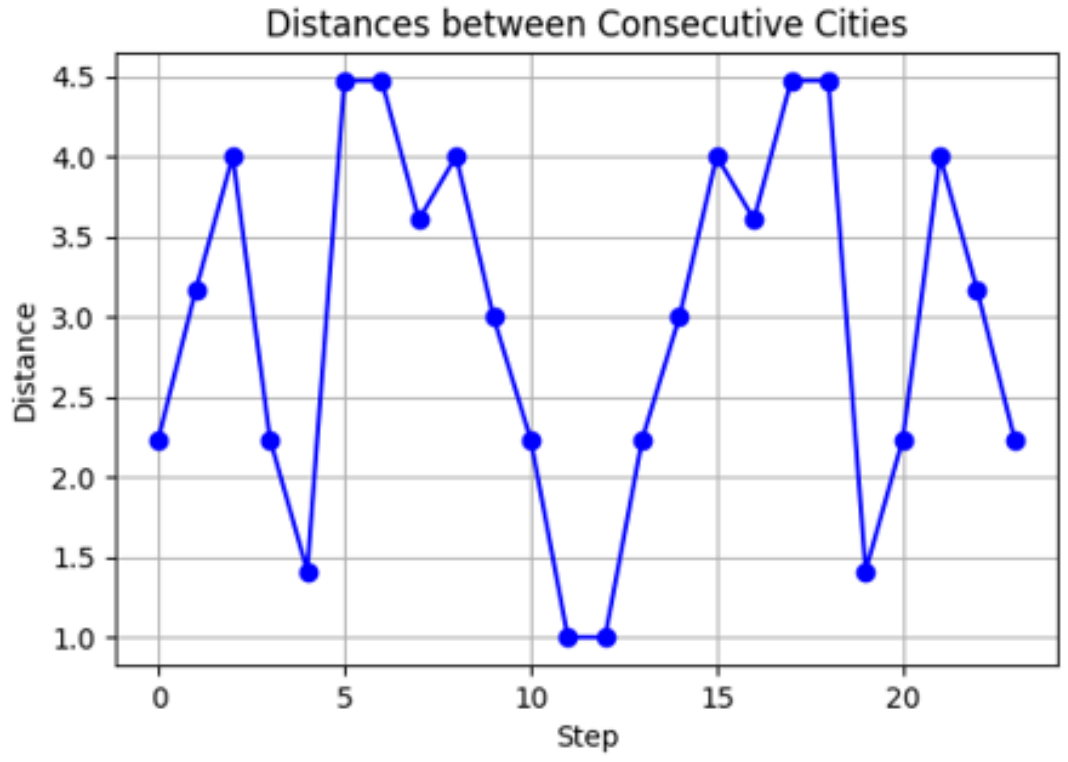

(c)

| | | | | |
|---|---|---|---|---|
| 6 | 7 | 23 | 17 | 12 |
| 5 | 15 | 25 | 16 | 4 |
| 18 | 24 | 13 | 2 | 8 |
| 22 | 10 | 1 | 11 | 21 |
| 14 | 9 | 3 | 19 | 20 |

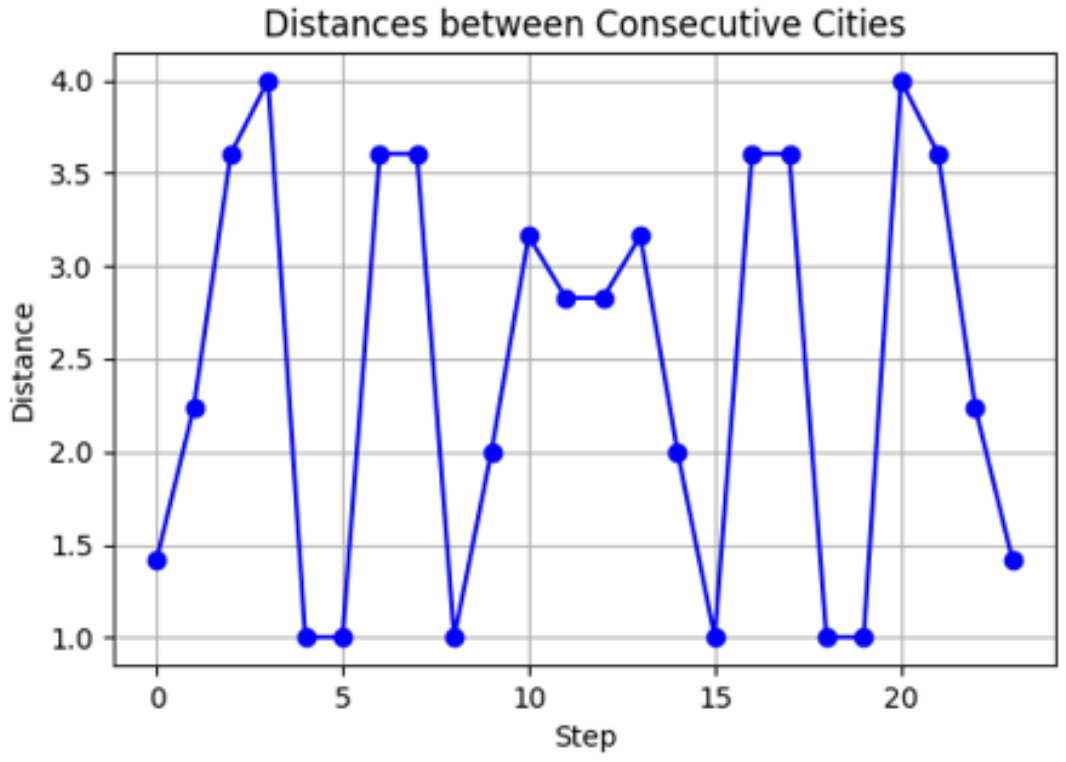

(d)

**Fig. 6:** Illustrating reflexive symmetry in distance patterns of the 5×5 associative magic path of the Traveler.

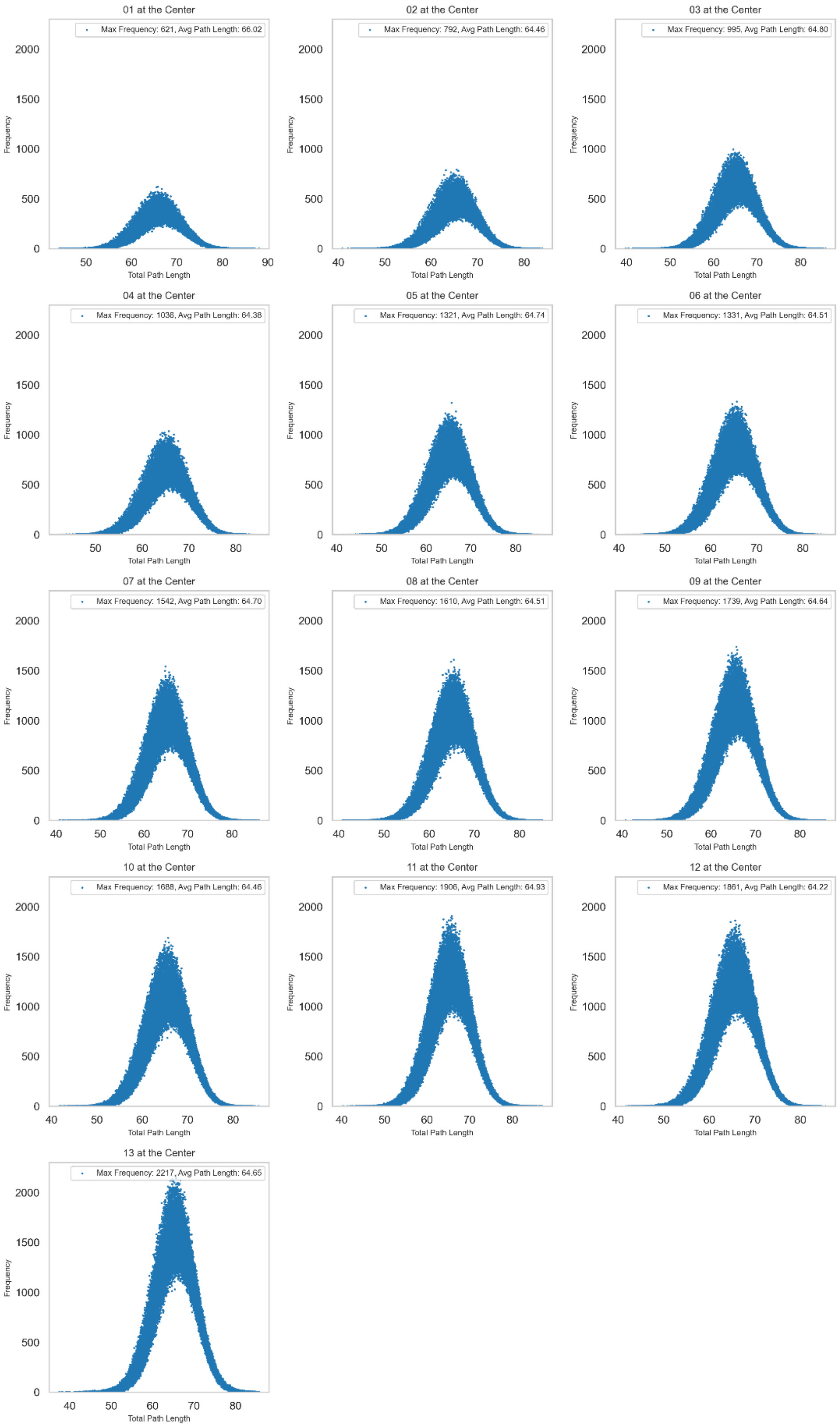

**Figure 7:** Distribution of total path lengths across different groups of order 5 magic squares, each with the same center number.

## 4. Concluding Remarks

In this study, an analysis of the topological distances of travelling salesman-like trajectories across the 3×3, 4×4, and 5×5 magic squares were investigated. In the smallest case, that of the unique irreducible 3×3 magic square, the lengths of trajectories exhibit clear reflective σ-symmetry.

A classification of 4×4 magic squares into the 12 groups suggested by Dudeney [3] reveals that a significant portion of the magic squares exhibit σ-symmetry in their distance patterns. Reflexive symmetry was observed in all magic squares categorized under Group 3 and Group 6, whereas other groups exhibited characteristics such as local symmetry, translational symmetry, and partial symmetry.

We investigated the total path length of all 5×5 magic squares, analyzing their frequencies and distributions. By classifying the squares into 25 groups, each with the same number at the center, we observed that groups with smaller center numbers tend to have lower maximum frequencies for total path length. Despite this variation, all 25 groups have an average total path length of approximately 64. Additionally, we identified a reflexive σ-symmetry in the distance patterns of both the associative 5×5 magic squares and 5,848 non-associative ones.

Future investigations are needed to explore the generalization of these findings to squares of higher orders. How can a traveling salesman find the shortest trajectory within 6×6 magic squares without considering all possible paths (approximately $1.7\times10^{19}$ paths)?

**Appendix: A listing of a python programme and algorithm used for order 4 magic squares**

```
import pandas as pd
import numpy as np
import matplotlib.pyplot as plt
import os

# Load the magic squares from the Excel file
df = pd.read_excel("magic_squares.xlsx", header=None)

# Define the city coordinates
cities = {
    1: (-1.5, 1.5),
    2: (-0.5, 1.5),
    3: (0.5, 1.5),
    4: (1.5, 1.5),
    5: (-1.5, 0.5),
    6: (-0.5, 0.5),
    7: (0.5, 0.5),
    8: (1.5, 0.5),
    9: (-1.5, -0.5),
    10: (-0.5, -0.5),
    11: (0.5, -0.5),
    12: (1.5, -0.5),
```

```
    13: (-1.5, -1.5),
    14: (-0.5, -1.5),
    15: (0.5, -1.5),
    16: (1.5, -1.5)
}

# Create a folder for saving the plots
folder_path = "Trajectory-Distance plots"
os.makedirs(folder_path, exist_ok=True)

# Calculate and plot the trajectory for each magic square
for i in range(len(df)):
    magic_square = df.iloc[i].values.tolist()

    # Find the maximum number in the magic square
    max_number = max(magic_square)

    # Create a dictionary to map numbers to city coordinates
    city_mapping = {magic_square[j]: cities[j+1] for j in range(max_number)}

    # Create the path based on the city coordinates
    path = [city_mapping[number] for number in range(1, max_number+1)]
    x = [coord[0] for coord in path]
    y = [coord[1] for coord in path]

    # Plot the trajectory with arrows
    for j in range(1, max_number):
        plt.arrow(x[j-1], y[j-1], x[j]-x[j-1], y[j]-y[j-1], head_width=0.15, head_length=0.15,
fc='black', ec='black')
    plt.title('City Trajectory')
    plt.xlabel('X')
    plt.ylabel('Y')
    plt.grid(True)
    plt.savefig(os.path.join(folder_path, f'trajectory_magic_square_{i+1}.png'), dpi=300)
    plt.close()

    # Calculate distances between consecutive cities
    distances = []
    for j in range(1, max_number):
        distance = np.sqrt((x[j] - x[j-1])**2 + (y[j] - y[j-1])**2)
        distances.append(distance)

    # Plot the distances
    plt.plot(distances, marker='o', linestyle='-', color='blue')
    plt.title('Distances between Consecutive Cities')
    plt.xlabel('Step')
```

```
    plt.ylabel('Distance')
    plt.grid(True)
    plt.savefig(os.path.join(folder_path, f'distances_magic_square_{i+1}.png'), dpi=300)
    plt.close()

    plt.show()
```

**Acknowledgments**
C.F.M. gratefully acknowledges funding from the Natural Sciences and Engineering Research Council of Canada (NSERC) and Mount Saint Vincent University.

**Conflict of Interest**
The author has no conflicts of interest to disclose.

**Supplementary Material**
The supplementary material accompanying this paper provides all the codes necessary to reproduce the results reported, a comprehensive collection of the trajectory and distance patterns for all magic squares of order 4, and the calculated total path lengths of all order 5 magic squares along with their distribution. It offers a detailed exploration of the symmetries, patterns, and characteristics exhibited by each magic square in terms of the traveler's path.

**References**

[1] W.S. Andrews, *Magic squares and cubes*, The Open Court Publishing Company, 2nd Ed., 1917.
[2] W.H. Benson, *New recreations with magic squares*, Dover Publications, 1st Ed., 1976.
[3] H.E. Dudeney, *Amusements in mathematics*, Dover Publications, 1st Ed., 1958.
[4] P. Fahimi, R. Javadi, *An introduction to magic squares and their physical applications*, ResearchGate (2016), 1-26.
[5] L. Sallows, The lost theorem, *Math. Intell*. **19** (1997) 51–54.
[6] A. Paz, Wonder cubes: theme and variations, *Math. Intell.* **44** (2022) 87–98.
[7] C. Boyer, Some notes on the magic squares of squares problem, *Math. Intell.* **27** (2005) 52–64.
[8] S. Cammann, The evolution of magic squares in China, *J. Am. Orient. Soc.* **80** (1960) 116–124.
[9] J. Sesiano, Magic squares in the tenth century: two Arabic treatises by Anṭaki and Buzjani, Springer, 1st Ed., 2017.
[10] L. Euler, Problema algebraicum ob affectiones prorsus singulares, *Opera Omnia* 1st Ser. **6** (1770) 287–315.
[11] W. Trump, *How many magic squares are there?*, (2024). http://www.trump.de/magic-squares/howmany.html.
[12] K. Pinn, C. Wieczerkowski, Number of magic squares from parallel tempering Monte Carlo, *Int. J. Mod. Phys. C.* **9** (1998) 541–546.
[13] G. Kato, S. Minato, Enumeration of associative magic squares of order 7, *J. Inf. Process* **28** (2020) 903–910.
[14] A. Kitajima, M. Kikuchi, Numerous but rare: an exploration of magic squares, *PLoS One* **10** (2015) e0125062.
[15] A. Rogers, P. Loly, The inertia tensor of a magic cube, *Am. J. Phys.* **72** (2004) 786–789.
[16] P. Loly, The invariance of the moment of inertia of magic squares, *Math. Gaz.* **88** (2004) 151–153.
[17] P. Fahimi, *Quasi-static levitation of magic squares*, Submitted, (2023).
[18] P. Fahimi, C.A. Toussi, W. Trump, J. Haddadnia, C.F. Matta, Striking patterns in natural magic squares' associated electrostatic potentials: Matrices of the 4th and 5th order, *Discrete Math.* **344** (2021) 112229.
[19] P. Fahimi, B. Jaleh, The electrostatic potential at the center of associative magic squares, *Int. J.*

*Phys. Sci*. **7** (2012) 24–30.
[20] A. Rogers, P. Loly, The electric multipole expansion for a magic cube, *Eur. J. Phys.* **26** (2005) 809-813.
[21] G. las Cuevas, T. Drescher, T. Netzer, Quantum magic squares: dilations and their limitations, *J. Math. Phys.* **61** (2020) 111704.
[22] G. las Cuevas, T. Netzer, I. Valentiner-Branth, Magic squares: Latin, semiclassical, and quantum, *J. Math. Phys.* **64** (2023) 22201.
[23] P. Fahimi, Binary color-coded magic squares: a study of uniqueness under rotation/reflection, PCA, and LDA analysis, *Discrete Math.* **347** (2024) 113708.
[24] L. Onsager, Crystal statistics. I. A two-dimensional model with an order-disorder transition, *Phys. Rev.* **65** (1944) 117-149.
[25] H. Heinz, *Order-4 magic squares*, (2008). http://recmath.org/Magic Squares/order4list.htm.